# ON NANO SEMI ALPHA OPEN SETS

QAYS HATEM IMRAN[1]

_________________________________________________

**Abstract.** *In this paper, we presented another concept of $N$-$O.S.$ called $NS_\alpha$-$O.S.$ and studied their fundamental properties in nano topological spaces. We also present $NS_\alpha$-interior and $NS_\alpha$-closure and study some of their fundamental properties.*
**Mathematics Subject Classification (2010):** *54A05, 54B05.*
**Keywords:** *$NS_\alpha$-$O.S.$, $NS_\alpha$-$C.S.$, $NS_\alpha$-interior and $NS_\alpha$-closure.*

## 1. INTRODUCTION

In 2000, G.B. Navalagi [1] presented the idea of semi-$\alpha$-open sets in topological spaces. N.M. Ali [2] introduced new types of weakly open sets in topological spaces. M.L. Thivagar and C. Richard [3] gave nano topological space (or simply N.T.S.) on a subset $\mathcal{M}$ of a universe which is defined regarding lower and upper approximations of $\mathcal{M}$. He studied about the weak forms of nano open sets (briefly $N$-$O.S.$), such as $N\alpha$-$O.S.$, $Ns$-$O.S.$, and $Np$-$O.S.$. The objective of this paper is to present the idea of $NS_\alpha$-$O.S.$ and study their fundamental properties in nano topological spaces. We also present $NS_\alpha$-interior and $NS_\alpha$-closure and obtain some of its properties.

## 2. PRELIMINARIES

Throughout this paper, $(\mathcal{U}, \tau_\mathcal{R}(\mathcal{M}))$ (or simply $\mathcal{U}$) always mean a nano topological space on which no separation axioms are expected unless generally specified. The complement of a $N$-$O.S.$ is called a nano closed set (briefly $N$-$C.S.$) in $(\mathcal{U}, \tau_\mathcal{R}(\mathcal{M}))$. For a set $\mathcal{C}$ in a nano topological space $(\mathcal{U}, \tau_\mathcal{R}(\mathcal{M}))$, $Ncl(\mathcal{C})$, $Nint(\mathcal{C})$ and $\mathcal{C}^c = \mathcal{U} - \mathcal{C}$ denote the nano closure of $\mathcal{C}$, the nano interior of $\mathcal{C}$ and the nano complement of $\mathcal{C}$ respectively.

**Definition 2.1 [3]:**

A subset $\mathcal{C}$ of an N.T.S. $(\mathcal{U}, \tau_\mathcal{R}(\mathcal{M}))$ is said to be:
(i) A nano pre-open set (briefly $Np$-$O.S.$) if $\mathcal{C} \subseteq Nint(Ncl(\mathcal{C}))$. The complement of a $Np$-$O.S.$ is called a nano pre-closed set (briefly $Np$-$C.S.$) in $(\mathcal{U}, \tau_\mathcal{R}(\mathcal{M}))$. The family of all $Np$-$O.S.$ (resp. $Np$-$C.S.$) of $\mathcal{U}$ is denoted by $NpO(\mathcal{U}, \mathcal{M})$ (resp. $NpC(\mathcal{U}, \mathcal{M})$).
(ii) A nano semi-open set (briefly $Ns$-$O.S.$) if $\mathcal{C} \subseteq Ncl(Nint(\mathcal{C}))$. The complement of a $Ns$-$O.S.$ is called a nano semi-closed set (briefly $Ns$-$C.S.$) in $(\mathcal{U}, \tau_\mathcal{R}(\mathcal{M}))$. The family of all $Ns$-$O.S.$ (resp. $Ns$-$C.S.$) of $\mathcal{U}$ is denoted by $NsO(\mathcal{U}, \mathcal{M})$ (resp. $NsC(\mathcal{U}, \mathcal{M})$).
(iii) A nano $\alpha$-open set (briefly $N\alpha$-$O.S.$) if $\mathcal{C} \subseteq Nint(Ncl(Nint(\mathcal{C})))$. The complement of a $N\alpha$-$O.S.$ is called a nano $\alpha$-closed set (briefly $N\alpha$-$C.S.$) in $(\mathcal{U}, \tau_\mathcal{R}(\mathcal{M}))$. The family of all $N\alpha$-$O.S.$ (resp. $N\alpha$-$C.S.$) of $\mathcal{U}$ is denoted by $N\alpha O(\mathcal{U}, \mathcal{M})$ (resp. $N\alpha C(\mathcal{U}, \mathcal{M})$).

_________________________________________________

[1] Muthanna University, College of Education for Pure Science, Department of Mathematics, Iraq.
E-mail: qays.imran@mu.edu.iq.

**Definition 2.2 [3]:**

(i) The $Np$-interior of a set $\mathcal{C}$ of a N.T.S. $(\mathcal{U}, \tau_{\mathcal{R}}(\mathcal{M}))$ is the union of all $Np$-O.S. contained in $\mathcal{C}$ and is denoted by $Npint(\mathcal{C})$.
(ii) The $Ns$-interior of a set $\mathcal{C}$ of a N.T.S. $(\mathcal{U}, \tau_{\mathcal{R}}(\mathcal{M}))$ is the union of all $Ns$-O.S. contained in $\mathcal{C}$ and is denoted by $Nsint(\mathcal{C})$.
(iii) The $N\alpha$-interior of a set $\mathcal{C}$ of a N.T.S. $(\mathcal{U}, \tau_{\mathcal{R}}(\mathcal{M}))$ is the union of all $N\alpha$-O.S. contained in $\mathcal{C}$ and is denoted by $N\alpha int(\mathcal{C})$.

**Definition 2.3 [3]:**

(i) The $Np$-closure of a set $\mathcal{C}$ of a N.T.S. $(\mathcal{U}, \tau_{\mathcal{R}}(\mathcal{M}))$ is the intersection of all $Np$-C.S. that contain $\mathcal{C}$ and is denoted by $Npcl(\mathcal{C})$.
(ii) The $Ns$-closure of a set $\mathcal{C}$ of a N.T.S. $(\mathcal{U}, \tau_{\mathcal{R}}(\mathcal{M}))$ is the intersection of all $Ns$-C.S. that contain $\mathcal{C}$ and is denoted by $Nscl(\mathcal{C})$.
(iii) The $N\alpha$-closure of a set $\mathcal{C}$ of a N.T.S. $(\mathcal{U}, \tau_{\mathcal{R}}(\mathcal{M}))$ is the intersection of all $N\alpha$-C.S. that contain $\mathcal{C}$ and is denoted by $N\alpha cl(\mathcal{C})$.

**Proposition 2.4 [3]:**

In a N.T.S. $(\mathcal{U}, \tau_{\mathcal{R}}(\mathcal{M}))$, then the following statements hold, and the equality of each statement are not true:
(i) Every $N$-O.S. (resp. $N$-C.S.) is a $N\alpha$-O.S. (resp. $N\alpha$-C.S.).
(ii) Every $N\alpha$-O.S. (resp. $N\alpha$-C.S.) is a $Ns$-O.S. (resp. $Ns$-C.S.).
(iii) Every $N\alpha$-O.S. (resp. $N\alpha$-C.S.) is a $Np$-O.S. (resp. $Np$-C.S.).

**Proposition 2.5 [3]:**

A subset $\mathcal{C}$ of a N.T.S. $(\mathcal{U}, \tau_{\mathcal{R}}(\mathcal{M}))$ is a $N\alpha$-O.S. iff $\mathcal{C}$ is a $Ns$-O.S. and $Np$-O.S..

**Lemma 2.6:**

(i) If $\mathcal{K}$ is a $N$-O.S., then $Nscl(\mathcal{K}) = Nint(Ncl(\mathcal{K}))$.
(ii) If $\mathcal{C}$ is a subset of a N.T.S. $(\mathcal{U}, \tau_{\mathcal{R}}(\mathcal{M}))$, then $Nsint(Ncl(\mathcal{C})) = Ncl(Nint(Ncl(\mathcal{C})))$.

### 3. NANO SEMI-$\alpha$-OPEN SETS

In this section, we present and study the $NS_\alpha$-O.S. and some of its properties.

**Definition 3.1:**

A subset $\mathcal{C}$ of a N.T.S. $(\mathcal{U}, \tau_{\mathcal{R}}(\mathcal{M}))$ is called nano semi-$\alpha$-open set (briefly $NS_\alpha$-O.S.) if there exists a $N\alpha$-O.S. $\mathcal{P}$ in $\mathcal{U}$ such that $\mathcal{P} \subseteq \mathcal{C} \subseteq Ncl(\mathcal{P})$ or equivalently if $\mathcal{C} \subseteq Ncl(N\alpha int(\mathcal{C}))$. The family of all $NS_\alpha$-O.S. of $\mathcal{U}$ is denoted by $NS_\alpha O(\mathcal{U}, \mathcal{M})$.

**Definition 3.2:**

The complement of $NS_\alpha$-O.S. is called a nano semi-$\alpha$-closed set (briefly $NS_\alpha$-C.S.). The family of all $NS_\alpha$-C.S. of $\mathcal{U}$ is denoted by $NS_\alpha C(\mathcal{U}, \mathcal{M})$.

**Example 3.3:**

Let $\mathcal{U} = \{p, q, r, s\}$ with $\mathcal{U}/\mathcal{R} = \{\{p\}, \{r\}, \{q, s\}\}$ and $\mathcal{M} = \{p, q\}$.
Let $\tau_{\mathcal{R}}(\mathcal{M}) = \{\phi, \{p\}, \{q, s\}, \{p, q, s\}, \mathcal{U}\}$ be a N.T.S.. The $N$-C.S. are $\mathcal{U}, \{q, r, s\}, \{p, r\}, \{r\}$ and $\phi$. The family of all $N\alpha$-O.S. of $\mathcal{U}$ is: $N\alpha O(\mathcal{U}, \mathcal{M}) = \{\phi, \{p\}, \{q, s\}, \{p, q, s\}, \mathcal{U}\}$.
The family of all $N\alpha$-C.S. of $\mathcal{U}$ is: $N\alpha C(\mathcal{U}, \mathcal{M}) = \{\mathcal{U}, \{q, r, s\}, \{p, r\}, \{r\}, \phi\}$.
The family of all $NS_\alpha$-O.S. of $\mathcal{U}$ is: $NS_\alpha O(\mathcal{U}, \mathcal{M}) = N\alpha O(\mathcal{U}, \mathcal{M}) \cup \{\{p, r\}, \{q, r, s\}\}$.
The family of all $NS_\alpha$-C.S. of $\mathcal{U}$ is: $NS_\alpha C(\mathcal{U}, \mathcal{M}) = N\alpha C(\mathcal{U}, \mathcal{M}) \cup \{\{q, s\}, \{p\}\}$.

**Remark 3.4:**

It is evident by definitions that in a N.T.S. $(\mathcal{U}, \tau_{\mathcal{R}}(\mathcal{M}))$, the following hold:
(i) Every $N$-O.S. (resp. $N$-C.S.) is a $NS_\alpha$-O.S. (resp. $NS_\alpha$-C.S.).
(ii) Every $N\alpha$-O.S. (resp. $N\alpha$-C.S.) is a $NS_\alpha$-O.S. (resp. $NS_\alpha$-C.S.).

The opposite of the above remark need not be true as appeared in the following example.

**Example 3.5:**

In example (3.3), the set $\{p, r\}$ is a $NS_\alpha$-O.S. but is not $N$-O.S. and not $N\alpha$-O.S.. The set $\{q, s\}$ is a $NS_\alpha$-C.S. but is not $N$-C.S. and not $N\alpha$-C.S..

**Remark 3.6:**

The concepts of $NS_\alpha$-O.S. and $Np$-O.S. are independent, as the following example shows.

**Example 3.7:**

In example (3.3), then the set $\{p, r\}$ is a $NS_\alpha$-O.S. but is not $Np$-O.S.. The set $\{p, r, s\}$ is a $Np$-O.S. but is not $NS_\alpha$-O.S..

**Remark 3.8:**

(i) If every $N$-O.S. is a $N$-C.S. and every nowhere nano dense set is $N$-C.S. in any N.T.S. $(\mathcal{U}, \tau_{\mathcal{R}}(\mathcal{M}))$, then every $NS_\alpha$-O.S. is a $N$-O.S..
(ii) If every $N$-O.S. is a $N$-C.S. in any N.T.S. $(\mathcal{U}, \tau_{\mathcal{R}}(\mathcal{M}))$, then every $NS_\alpha$-O.S. is a $N\alpha$-O.S..

**Remark 3.9:**

(i) It is clear that every $Ns$-O.S. and $Np$-O.S. of any N.T.S. $(\mathcal{U}, \tau_{\mathcal{R}}(\mathcal{M}))$ is a $NS_\alpha$-O.S. (by proposition (2.5) and remark (3.4) (ii)).
(ii) A $NS_\alpha$-O.S. in any N.T.S. $(\mathcal{U}, \tau_{\mathcal{R}}(\mathcal{M}))$ is a $Np$-O.S. if every $N$-O.S. of $\mathcal{U}$ is a $N$-C.S. (from proposition (2.4) (iii) and remark (3.8) (ii)).

**Theorem 3.10:**

For any subset $\mathcal{C}$ of a N.T.S. $(\mathcal{U}, \tau_{\mathcal{R}}(\mathcal{M}))$, $\mathcal{C} \in N\alpha O(\mathcal{U}, \mathcal{M})$ iff there exists a $N$-O.S. $\mathcal{K}$ such that $\mathcal{K} \subseteq \mathcal{C} \subseteq Nint(Ncl(\mathcal{K}))$.

*Proof:* Let $\mathcal{C}$ be a $N\alpha$-O.S.. Hence $\mathcal{C} \subseteq Nint(Ncl(Nint(\mathcal{C})))$, so let $\mathcal{K} = Nint(\mathcal{C})$, we get $Nint(\mathcal{C}) \subseteq \mathcal{C} \subseteq Nint(Ncl(Nint(\mathcal{C})))$. Then there exists a $N$-O.S. $Nint(\mathcal{C})$ such that $\mathcal{K} \subseteq \mathcal{C} \subseteq Nint(Ncl(\mathcal{K}))$, where $\mathcal{K} = Nint(\mathcal{C})$.

Conversely, suppose that there is a $N$-O.S. $\mathcal{K}$ such that $\mathcal{K} \subseteq \mathcal{C} \subseteq Nint(Ncl(\mathcal{K}))$.
To prove $\mathcal{C} \in N\alpha O(\mathcal{U}, \mathcal{M})$.
$\mathcal{K} \subseteq Nint(\mathcal{C})$ (since $Nint(\mathcal{C})$ is the largest $N$-O.S. contained in $\mathcal{C}$).
Hence $Ncl(\mathcal{K}) \subseteq Nint(Ncl(\mathcal{C}))$, then $Nint(Ncl(\mathcal{K})) \subseteq Nint(Ncl(Nint(\mathcal{C})))$.
But $\mathcal{K} \subseteq \mathcal{C} \subseteq Nint(Ncl(\mathcal{K}))$ (by hypothesis). Then $\mathcal{C} \subseteq Nint(Ncl(Nint(\mathcal{C})))$.
Therefore, $\mathcal{C} \in N\alpha O(\mathcal{U}, \mathcal{M})$.

**Theorem 3.11:**

For any subset $\mathcal{C}$ of a N.T.S. $(\mathcal{U}, \tau_\mathcal{R}(\mathcal{M}))$. The following properties are equivalent:
(i) $\mathcal{C} \in NS_\alpha O(\mathcal{U}, \mathcal{M})$.
(ii) There exists a $N$-O.S. say $\mathcal{K}$ such that $\mathcal{K} \subseteq \mathcal{C} \subseteq Ncl(Nint(Ncl(\mathcal{K})))$.
(iii) $\mathcal{C} \subseteq Ncl(Nint(Ncl(Nint(\mathcal{C}))))$.

*Proof:*
$(i) \Rightarrow (ii)$ Let $\mathcal{C} \in NS_\alpha O(\mathcal{U}, \mathcal{M})$. Then there exists $\mathcal{P} \in N\alpha O(\mathcal{U}, \mathcal{M})$, such that $\mathcal{P} \subseteq \mathcal{C} \subseteq Ncl(\mathcal{P})$. Hence there exists $\mathcal{K}$ $N$-O.S. such that $\mathcal{K} \subseteq \mathcal{P} \subseteq Nint(Ncl(\mathcal{K}))$ (by theorem (3.10)). Therefore, $Ncl(\mathcal{K}) \subseteq Ncl(\mathcal{P}) \subseteq Ncl(Nint(Ncl(\mathcal{K})))$, implies that $Ncl(\mathcal{P}) \subseteq Ncl(Nint(Ncl(\mathcal{K})))$. Then $\mathcal{K} \subseteq \mathcal{P} \subseteq \mathcal{C} \subseteq Ncl(\mathcal{P}) \subseteq Ncl(Nint(Ncl(\mathcal{K})))$. Therefore, $\mathcal{K} \subseteq \mathcal{C} \subseteq Ncl(Nint(Ncl(\mathcal{K})))$, for some $\mathcal{K}$ $N$-O.S..
$(ii) \Rightarrow (iii)$ Suppose that there exists a $N$-O.S. $\mathcal{K}$ such that $\mathcal{K} \subseteq \mathcal{C} \subseteq Ncl(Nint(Ncl(\mathcal{K})))$. We know that $Nint(\mathcal{C}) \subseteq \mathcal{C}$. On the other hand, $\mathcal{K} \subseteq Nint(\mathcal{C})$ (since $Nint(\mathcal{C})$ is the largest $N$-O.S. contained in $\mathcal{C}$). Hence $Ncl(\mathcal{K}) \subseteq Ncl(Nint(\mathcal{C}))$, then $Nint(Ncl(\mathcal{K})) \subseteq Nint(Ncl(Nint(\mathcal{C})))$, therefore $Ncl(Nint(Ncl(\mathcal{K}))) \subseteq Ncl(Nint(Ncl(Nint(\mathcal{C}))))$.
But $\mathcal{C} \subseteq Ncl(Nint(Ncl(\mathcal{K})))$ (by hypothesis).
Hence $\mathcal{C} \subseteq Ncl(Nint(Ncl(\mathcal{K}))) \subseteq Ncl(Nint(Ncl(Nint(\mathcal{C}))))$,
then $\mathcal{C} \subseteq Ncl(Nint(Ncl(Nint(\mathcal{C}))))$.
$(iii) \Rightarrow (i)$ Let $\mathcal{C} \subseteq Ncl(Nint(Ncl(Nint(\mathcal{C}))))$. To prove $\mathcal{C} \in NS_\alpha O(\mathcal{U}, \mathcal{M})$.
Let $\mathcal{P} = Nint(\mathcal{C})$; we know that $Nint(\mathcal{C}) \subseteq \mathcal{C}$. To prove $\mathcal{C} \subseteq Ncl(Nint(\mathcal{C}))$.
Since $Nint(Ncl(Nint(\mathcal{C}))) \subseteq Ncl(Nint(\mathcal{C}))$.
Hence, $Ncl(Nint(Ncl(Nint(\mathcal{C})))) \subseteq Ncl(Ncl(Nint(\mathcal{C}))) = Ncl(Nint(\mathcal{C}))$.
But $\mathcal{C} \subseteq Ncl(Nint(Ncl(Nint(\mathcal{C}))))$ (by hypothesis).
Hence, $\mathcal{C} \subseteq Ncl(Nint(Ncl(Nint(\mathcal{C})))) \subseteq Ncl(Nint(\mathcal{C})) \Rightarrow \mathcal{C} \subseteq Ncl(Nint(\mathcal{C}))$.
Hence, there exists a $N$-O.S. say $\mathcal{P}$, such that $\mathcal{P} \subseteq \mathcal{C} \subseteq Ncl(\mathcal{P})$.
On the other hand, $\mathcal{P}$ is a $N\alpha$-O.S. (since $\mathcal{P}$ is a $N$-O.S.). Hence $\mathcal{C} \in NS_\alpha O(\mathcal{U}, \mathcal{M})$.

**Corollary 3.12:**

For any subset $\mathcal{C}$ of a N.T.S. $(\mathcal{U}, \tau_\mathcal{R}(\mathcal{M}))$, the following properties are equivalent:
(i) $\mathcal{C} \in NS_\alpha C(\mathcal{U}, \mathcal{M})$.
(ii) There exists a $N$-C.S. $\mathcal{F}$ such that $Nint(Ncl(Nint(\mathcal{F}))) \subseteq \mathcal{C} \subseteq \mathcal{F}$.
(iii) $Nint(Ncl(Nint(Ncl(\mathcal{C})))) \subseteq \mathcal{C}$.

*Proof:*
$(i) \Rightarrow (ii)$ Let $\mathcal{C} \in NS_\alpha C(\mathcal{U}, \mathcal{M})$, then $\mathcal{C}^c \in NS_\alpha O(\mathcal{U}, \mathcal{M})$. Hence there is $\mathcal{K}$ $N$-O.S. such that $\mathcal{K} \subseteq \mathcal{C}^c \subseteq Ncl(Nint(Ncl(\mathcal{K})))$ (by theorem (3.11)). Hence $(Ncl(Nint(Ncl(\mathcal{K}))))^c \subseteq \mathcal{C}^{cc} \subseteq \mathcal{K}^c$, i.e., $Nint(Ncl(Nint(\mathcal{K}^c))) \subseteq \mathcal{C} \subseteq \mathcal{K}^c$. Let $\mathcal{K}^c = \mathcal{F}$, where $\mathcal{F}$ is a $N$-C.S. in $\mathcal{U}$. Then $Nint(Ncl(Nint(\mathcal{F}))) \subseteq \mathcal{C} \subseteq \mathcal{F}$, for some $\mathcal{F}$ $N$-C.S..
$(ii) \Rightarrow (iii)$ Suppose that there exists $\mathcal{F}$ $N$-C.S. such that $Nint\big(Ncl\big(Nint(\mathcal{F})\big)\big) \subseteq \mathcal{C} \subseteq \mathcal{F}$, but $Ncl(\mathcal{C})$ is the smallest $N$-C.S. containing $\mathcal{C}$. Then $Ncl(\mathcal{C}) \subseteq \mathcal{F}$, and therefore:

$Nint(Ncl(C)) \subseteq Nint(\mathcal{F}) \Rightarrow Ncl\big(Nint(Ncl(C))\big) \subseteq Ncl(Nint(\mathcal{F})) \Rightarrow$
$Nint(Ncl(Nint(Ncl(C)))) \subseteq Nint(Ncl(Nint(\mathcal{F}))) \subseteq C \Rightarrow$
$Nint(Ncl(Nint(Ncl(C)))) \subseteq C$.

$(iii) \Rightarrow (i)$ Let $Nint(Ncl(Nint(Ncl(C)))) \subseteq C$. To prove $C \in NS_\alpha C(\mathcal{U}, \mathcal{M})$,
i.e., to prove $C^c \in NS_\alpha O(\mathcal{U}, \mathcal{M})$.
Then $C^c \subseteq (Nint(Ncl(Nint(Ncl(C)))))^c = Ncl(Nint(Ncl(Nint(C^c))))$, but
$(Nint(Ncl(Nint(Ncl(C)))))^c = Ncl(Nint(Ncl(Nint(C^c))))$.
Hence $C^c \subseteq Ncl(Nint(Ncl(Nint(C^c))))$, and therefore $C^c \in NS_\alpha O(\mathcal{U}, \mathcal{M})$,
i.e., $C \in NS_\alpha C(\mathcal{U}, \mathcal{M})$.

**Proposition 3.13:**

The union of any family of $N\alpha$-O.S. is a $N\alpha$-O.S..

*Proof:* Let $\{C_i\}_{i \in \Lambda}$ be a family of $N\alpha$-O.S. of $\mathcal{U}$. To prove $\bigcup_{i \in \Lambda} C_i$ is a $N\alpha$-O.S.,
i.e., $\bigcup_{i \in \Lambda} C_i \subseteq Nint(Ncl(Nint(\bigcup_{i \in \Lambda} C_i)))$. Then $C_i \subseteq Nint(Ncl(Nint(C_i)))$, $\forall i \in \Lambda$.
Since $\bigcup_{i \in \Lambda} Nint(C_i) \subseteq Nint(\bigcup_{i \in \Lambda} C_i)$ and $\bigcup_{i \in \Lambda} Ncl(C_i) \subseteq Ncl(\bigcup_{i \in \Lambda} C_i)$ hold for any nano topology. We have $\bigcup_{i \in \Lambda} C_i \subseteq \bigcup_{i \in \Lambda} Nint(Ncl(Nint(C_i)))$
$\subseteq Nint(\bigcup_{i \in \Lambda} Ncl(Nint(C_i)))$
$\subseteq Nint(Ncl(\bigcup_{i \in \Lambda}(Nint(C_i)))$
$\subseteq Nint(Ncl(Nint(\bigcup_{i \in \Lambda} C_i)))$.
Hence $\bigcup_{i \in \Lambda} C_i$ is a $N\alpha$-O.S..

**Theorem 3.14:**

The union of any family of $NS_\alpha$-O.S. is a $NS_\alpha$-O.S..

*Proof:* Let $\{C_i\}_{i \in \Lambda}$ be a family of $NS_\alpha$-O.S.. To prove $\bigcup_{i \in \Lambda} C_i$ is a $NS_\alpha$-O.S..
Since $C_i \in NS_\alpha O(\mathcal{U}, \mathcal{M})$. Then there is a $N\alpha$-O.S. $\mathcal{D}_i$ such that $\mathcal{D}_i \subseteq C_i \subseteq Ncl(\mathcal{D}_i)$, $\forall i \in \Lambda$.
Hence $\bigcup_{i \in \Lambda} \mathcal{D}_i \subseteq \bigcup_{i \in \Lambda} C_i \subseteq \bigcup_{i \in \Lambda} Ncl(\mathcal{D}_i) \subseteq Ncl(\bigcup_{i \in \Lambda} \mathcal{D}_i)$.
But $\bigcup_{i \in \Lambda} \mathcal{D}_i \in N\alpha O(\mathcal{U}, \mathcal{M})$ (by proposition (3.13)). Hence $\bigcup_{i \in \Lambda} C_i \in NS_\alpha O(\mathcal{U}, \mathcal{M})$.

**Corollary 3.15:**

The intersection of any family of $NS_\alpha$-C.S. is a $NS_\alpha$-C.S..

*Proof:* This follows directly from the theorem (3.14).

**Remark 3.16:**

The intersection of any two $NS_\alpha$-O.S. is not necessary $NS_\alpha$-O.S. as in the following example.

**Example 3.17:**

In example (3.3), $\{p, r\}$ and $\{q, r, s\}$ are two $NS_\alpha$-O.S., but $\{p, r\} \cap \{q, r, s\} = \{r\}$ is not $NS_\alpha$-O.S..

**Remark 3.18:**

The following diagram shows the relations among the different types of weakly $N$-O.S. that were studied in this section:

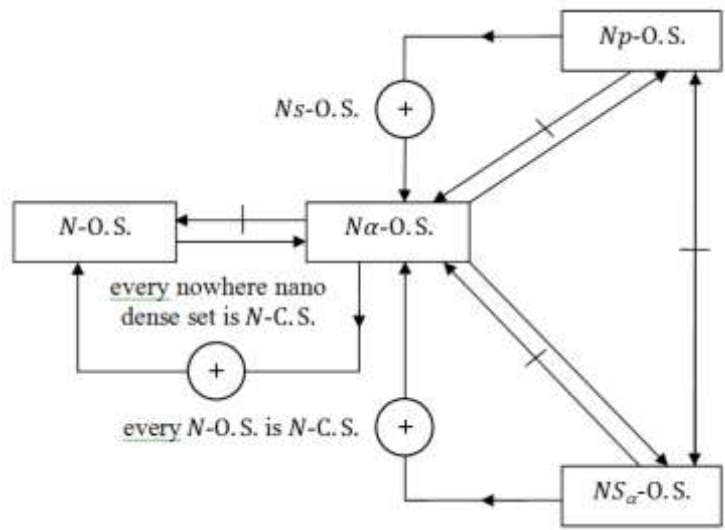

## 4. NANO SEMI-$\alpha$-INTERIOR AND NANO SEMI-$\alpha$-CLOSURE

We present $NS_\alpha$-interior and $NS_\alpha$-closure and obtain some of its properties in this section.

**Definition 4.1:**

The union of all $NS_\alpha$-O.S. in a N.T.S. $(\mathcal{U}, \tau_\mathcal{R}(\mathcal{M}))$ contained in $\mathcal{C}$ is called $NS_\alpha$-interior of $\mathcal{C}$ and is denoted by $NS_\alpha int(\mathcal{C})$, $NS_\alpha int(\mathcal{C}) = \bigcup\{\mathcal{D}: \mathcal{D} \subseteq \mathcal{C}, \mathcal{D}$ is a $NS_\alpha$-O.S.$\}$.

**Definition 4.2:**

The intersection of all $NS_\alpha$-C.S. in a N.T.S. $(\mathcal{U}, \tau_\mathcal{R}(\mathcal{M}))$ containing $\mathcal{C}$ is called $NS_\alpha$-closure of $\mathcal{C}$ and is denoted by $NS_\alpha cl(\mathcal{C})$, $NS_\alpha cl(\mathcal{C}) = \bigcap\{\mathcal{D}: \mathcal{C} \subseteq \mathcal{D}, \mathcal{D}$ is a $NS_\alpha$-C.S.$\}$.

**Proposition 4.3:**

Let $\mathcal{C}$ be any set in a N.T.S. $(\mathcal{U}, \tau_\mathcal{R}(\mathcal{M}))$, the following properties are true:
(i) $NS_\alpha int(\mathcal{C}) = \mathcal{C}$ iff $\mathcal{C}$ is a $NS_\alpha$-O.S..
(ii) $NS_\alpha cl(\mathcal{C}) = \mathcal{C}$ iff $\mathcal{C}$ is a $NS_\alpha$-C.S..
(iii) $NS_\alpha int(\mathcal{C})$ is the largest $NS_\alpha$-O.S. contained in $\mathcal{C}$.
(iv) $NS_\alpha cl(\mathcal{C})$ is the smallest $NS_\alpha$-C.S. containing $\mathcal{C}$.
*Proof:* (i), (ii), (iii) and (iv) are obvious.

**Proposition 4.4:**

Let $\mathcal{C}$ be any set in a N.T.S. $(\mathcal{U}, \tau_\mathcal{R}(\mathcal{M}))$, the following properties are true:
(i) $NS_\alpha int(\mathcal{U} - \mathcal{C}) = \mathcal{U} - (NS_\alpha cl(\mathcal{C}))$,
(ii) $NS_\alpha cl(\mathcal{U} - \mathcal{C}) = \mathcal{U} - (NS_\alpha int(\mathcal{C}))$.
*Proof:* (i) By definition, $NS_\alpha cl(\mathcal{C}) = \bigcap\{\mathcal{D}: \mathcal{C} \subseteq \mathcal{D}, \mathcal{D}$ is a $NS_\alpha$-C.S.$\}$
$\mathcal{U} - (NS_\alpha cl(\mathcal{C})) = \mathcal{U} - \bigcap\{\mathcal{D}: \mathcal{C} \subseteq \mathcal{D}, \mathcal{D}$ is a $NS_\alpha$-C.S.$\}$
$\qquad = \bigcup\{\mathcal{U} - \mathcal{D}: \mathcal{C} \subseteq \mathcal{D}, \mathcal{D}$ is a $NS_\alpha$-C.S.$\}$
$\qquad = \bigcup\{\mathcal{H}: \mathcal{H} \subseteq \mathcal{U} - \mathcal{C}, \mathcal{H}$ is a $NS_\alpha$-O.S.$\}$
$\qquad = NS_\alpha int(\mathcal{U} - \mathcal{C})$.

(ii) The proof is similar to (i).

**Theorem 4.5:**

Let $\mathcal{C}$ and $\mathcal{D}$ be two sets in a N.T.S. $(\mathcal{U}, \tau_\mathcal{R}(\mathcal{M}))$. The following properties hold:
(i) $NS_\alpha int(\phi) = \phi, NS_\alpha int(\mathcal{U}) = \mathcal{U}$.
(ii) $NS_\alpha int(\mathcal{C}) \subseteq \mathcal{C}$.
(iii) $\mathcal{C} \subseteq \mathcal{D} \Longrightarrow NS_\alpha int(\mathcal{C}) \subseteq NS_\alpha int(\mathcal{D})$.
(iv) $NS_\alpha int(\mathcal{C} \cap \mathcal{D}) \subseteq NS_\alpha int(\mathcal{C}) \cap NS_\alpha int(\mathcal{D})$.
(v) $NS_\alpha int(\mathcal{C}) \cup NS_\alpha int(\mathcal{D}) \subseteq NS_\alpha int(\mathcal{C} \cup \mathcal{D})$.
(vi) $NS_\alpha int(NS_\alpha int(\mathcal{C})) = NS_\alpha int(\mathcal{C})$.

*Proof:* (i), (ii), (iii), (iv), (v) and (vi) are obvious.

The equality in (iv) and (v) is not true in general, as the following example shows:

**Example 4.6:**

Let $\mathcal{U} = \{p, q, r, s\}$ with $\mathcal{U}/\mathcal{R} = \{\{q\}, \{r\}, \{p, s\}\}$ and $\mathcal{M} = \{p, r\}$.
Let $\tau_\mathcal{R}(\mathcal{M}) = \{\phi, \{r\}, \{p, s\}, \{p, r, s\}, \mathcal{U}\}$ be a N.T.S.. The $N$-C.S. are $\mathcal{U}, \{p, q, s\}, \{q, r\}, \{q\}$ and $\phi$. The family of all $N\alpha$-O.S. of $\mathcal{U}$ is: $N\alpha O(\mathcal{U}, \mathcal{M}) = \{\phi, \{r\}, \{p, s\}, \{p, r, s\}, \mathcal{U}\}$.
The family of all $NS_\alpha$-O.S. of $\mathcal{U}$ is: $NS_\alpha O(\mathcal{U}, \mathcal{M}) = N\alpha O(\mathcal{U}, \mathcal{M}) \cup \{\{q, r\}, \{p, q, s\}\}$.
Let $\mathcal{C} = \{q, r\}, \mathcal{D} = \{p, q, s\}$. Then $NS_\alpha int(\mathcal{C}) = \{q, r\}, NS_\alpha int(\mathcal{D}) = \{p, q, s\}, \mathcal{C} \cap \mathcal{D} = \{q\}$, $NS_\alpha int(\mathcal{C} \cap \mathcal{D}) = \phi$ and $NS_\alpha int(\mathcal{C}) \cap NS_\alpha int(\mathcal{D}) = \{q\}$.
It is clear that $NS_\alpha int(\mathcal{C}) \cap NS_\alpha int(\mathcal{D}) \nsubseteq NS_\alpha int(\mathcal{C} \cap \mathcal{D})$.
Let $\mathcal{C} = \{p, s\}, \mathcal{D} = \{q, s\}$. Then $NS_\alpha int(\mathcal{C}) = \{p, s\}, NS_\alpha int(\mathcal{D}) = \phi, \mathcal{C} \cup \mathcal{D} = \{p, q, s\}$, $NS_\alpha int(\mathcal{C} \cup \mathcal{D}) = \{p, q, s\}$ and $NS_\alpha int(\mathcal{C}) \cup NS_\alpha int(\mathcal{D}) = \{p, s\}$.
It is clear that $NS_\alpha int(\mathcal{C} \cup \mathcal{D}) \nsubseteq NS_\alpha int(\mathcal{C}) \cup NS_\alpha int(\mathcal{D})$.

**Theorem 4.7:**

Let $\mathcal{C}$ and $\mathcal{D}$ be two sets in a N.T.S. $(\mathcal{U}, \tau_\mathcal{R}(\mathcal{M}))$. The following properties hold:
(i) $NS_\alpha cl(\phi) = \phi, NS_\alpha cl(\mathcal{U}) = \mathcal{U}$.
(ii) $\mathcal{C} \subseteq NS_\alpha cl(\mathcal{C})$.
(iii) $\mathcal{C} \subseteq \mathcal{D} \Longrightarrow NS_\alpha cl(\mathcal{C}) \subseteq NS_\alpha cl(\mathcal{D})$.
(iv) $NS_\alpha cl(\mathcal{C} \cap \mathcal{D}) \subseteq NS_\alpha cl(\mathcal{C}) \cap NS_\alpha cl(\mathcal{D})$.
(v) $NS_\alpha cl(\mathcal{C}) \cup NS_\alpha cl(\mathcal{D}) \subseteq NS_\alpha cl(\mathcal{C} \cup \mathcal{D})$.
(vi) $NS_\alpha cl(NS_\alpha cl(\mathcal{C})) = NS_\alpha cl(\mathcal{C})$.

*Proof:* (i) and (ii) are evident.
(iii) By part (ii), $\mathcal{D} \subseteq NS_\alpha cl(\mathcal{D})$. Since $\mathcal{C} \subseteq \mathcal{D}$, we have $\mathcal{C} \subseteq NS_\alpha cl(\mathcal{D})$. But $NS_\alpha cl(\mathcal{D})$ is a $NS_\alpha$-C.S.. Thus $NS_\alpha$-$cl(\mathcal{D})$ is a $NS_\alpha$-C.S. containing $\mathcal{C}$. Since $NS_\alpha cl(\mathcal{C})$ is the smallest $NS_\alpha$-C.S. containing $\mathcal{C}$, we have $NS_\alpha cl(\mathcal{C}) \subseteq NS_\alpha cl(\mathcal{D})$. Hence, $\mathcal{C} \subseteq \mathcal{D} \Longrightarrow NS_\alpha cl(\mathcal{C}) \subseteq NS_\alpha cl(\mathcal{D})$.
(iv) We know that $\mathcal{C} \cap \mathcal{D} \subseteq \mathcal{C}$ and $\mathcal{C} \cap \mathcal{D} \subseteq \mathcal{D}$. Therefore, by part (iii), $NS_\alpha cl(\mathcal{C} \cap \mathcal{D}) \subseteq NS_\alpha cl(\mathcal{C})$ and $S_\alpha cl(\mathcal{C} \cap \mathcal{D}) \subseteq NS_\alpha cl(\mathcal{D})$. Hence $NS_\alpha cl(\mathcal{C} \cap \mathcal{D}) \subseteq NS_\alpha cl(\mathcal{C}) \cap NS_\alpha cl(\mathcal{D})$.
(v) Since $\mathcal{C} \subseteq \mathcal{C} \cup \mathcal{D}$ and $\mathcal{D} \subseteq \mathcal{C} \cup \mathcal{D}$, it follows from part (iii) that $NS_\alpha cl(\mathcal{C}) \subseteq NS_\alpha cl(\mathcal{C} \cup \mathcal{D})$ and $NS_\alpha cl(\mathcal{D}) \subseteq NS_\alpha cl(\mathcal{C} \cup \mathcal{D})$. Hence $NS_\alpha cl(\mathcal{C}) \cup NS_\alpha cl(\mathcal{D}) \subseteq NS_\alpha cl(\mathcal{C} \cup \mathcal{D})$.
(vi) Since $NS_\alpha cl(\mathcal{C})$ is a $NS_\alpha$-C.S., we have by proposition (4.3) part (ii), $S_\alpha cl(NS_\alpha cl(\mathcal{C})) = NS_\alpha cl(\mathcal{C})$.

The equality in (iv) and (v) is not true in general, as the following example shows:

**Example 4.8:**

In example (4.6), the family of all $N\alpha$-C.S. of $\mathcal{U}$ is: $N\alpha C(\mathcal{U},\mathcal{M}) = \{\mathcal{U},\{p,q,s\},\{q,r\},\{q\},\phi\}$. The family of all $NS_\alpha$-C.S. of $\mathcal{U}$ is: $NS_\alpha C(\mathcal{U},\mathcal{M}) = N\alpha C(\mathcal{U},\mathcal{M}) \cup \{\{p,s\},\{r\}\}$. Let $\mathcal{C} = \{p,r\}, \mathcal{D} = \{q,r\}$. Then $NS_\alpha cl(\mathcal{C}) = \mathcal{U}, NS_\alpha cl(\mathcal{D}) = \{q,r\}, \mathcal{C} \cap \mathcal{D} = \{r\}, NS_\alpha cl(\mathcal{C} \cap \mathcal{D}) = \{r\}$ and $NS_\alpha cl(\mathcal{C}) \cap NS_\alpha cl(\mathcal{D}) = \{q,r\}$.
It is clear that $NS_\alpha cl(\mathcal{C}) \cap NS_\alpha cl(\mathcal{D}) \nsubseteq NS_\alpha cl(\mathcal{C} \cap \mathcal{D})$.
Let $\mathcal{C} = \{p,s\}, \mathcal{D} = \{r\}$. Then $NS_\alpha cl(\mathcal{C}) = \{p,s\}, NS_\alpha cl(\mathcal{D}) = \{r\}, \mathcal{C} \cup \mathcal{D} = \{p,r,s\}$, $NS_\alpha cl(\mathcal{C} \cup \mathcal{D}) = \mathcal{U}$ and $NS_\alpha cl(\mathcal{C}) \cup NS_\alpha cl(\mathcal{D}) = \{p,r,s\}$.
It is clear that $NS_\alpha cl(\mathcal{C} \cup \mathcal{D}) \nsubseteq NS_\alpha cl(\mathcal{C}) \cup NS_\alpha cl(\mathcal{D})$.

**Proposition 4.9:**

For any subset $\mathcal{C}$ of a N.T.S. $(\mathcal{U}, \tau_\mathcal{R}(\mathcal{M}))$, then:
(i) $Nint(\mathcal{C}) \subseteq N\alpha int(\mathcal{C}) \subseteq NS_\alpha int(\mathcal{C}) \subseteq NS_\alpha cl(\mathcal{C}) \subseteq N\alpha cl(\mathcal{C}) \subseteq Ncl(\mathcal{C})$.
(ii) $Nint(NS_\alpha int(\mathcal{C})) = NS_\alpha int(Nint(\mathcal{C})) = Nint(\mathcal{C})$.
(iii) $N\alpha int(NS_\alpha int(\mathcal{C})) = NS_\alpha int(N\alpha int(\mathcal{C})) = N\alpha int(\mathcal{C})$.
(iv) $Ncl(NS_\alpha cl(\mathcal{C})) = NS_\alpha cl(Ncl(\mathcal{C})) = Ncl(\mathcal{C})$.
(v) $N\alpha cl(NS_\alpha cl(\mathcal{C})) = NS_\alpha cl(N\alpha cl(\mathcal{C})) = N\alpha cl(\mathcal{C})$.
(vi) $NS_\alpha cl(\mathcal{C}) = \mathcal{C} \cup Nint(Ncl(Nint(Ncl(\mathcal{C}))))$.
(vii) $NS_\alpha int(\mathcal{C}) = \mathcal{C} \cap Ncl(Nint(Ncl(Nint(\mathcal{C}))))$.
(viii) $Nint(Ncl(\mathcal{C})) \subseteq NS_\alpha int(NS_\alpha cl(\mathcal{C}))$.

*Proof:* We shall prove only (ii), (iii), (iv), (vii) and (viii).
(ii) To prove $Nint(NS_\alpha int(\mathcal{C})) = NS_\alpha int(Nint(\mathcal{C})) = Nint(\mathcal{C})$.
Since $Nint(\mathcal{C})$ is a $N$-O.S., then $Nint(\mathcal{C})$ is a $NS_\alpha$-O.S..
Hence $Nint(\mathcal{C}) = NS_\alpha int(Nint(\mathcal{C}))$ (by proposition (4.3)). Therefore:
$Nint(\mathcal{C}) = NS_\alpha int(Nint(\mathcal{C}))$ ................................................................................. (1)
Since $Nint(\mathcal{C}) \subseteq NS_\alpha int(\mathcal{C}) \Rightarrow Nint(Nint(\mathcal{C})) \subseteq Nint(NS_\alpha int(\mathcal{C})) \Rightarrow Nint(\mathcal{C}) \subseteq Nint(NS_\alpha int(\mathcal{C}))$.
Also, $NS_\alpha int(\mathcal{C}) \subseteq \mathcal{C} \Rightarrow Nint(NS_\alpha int(\mathcal{C})) \subseteq Nint(\mathcal{C})$. Hence:
$Nint(\mathcal{C}) = Nint(NS_\alpha int(\mathcal{C}))$ ................................................................................. (2)
Therefore by (1) and (2), we get $Nint(NS_\alpha int(\mathcal{C})) = NS_\alpha int(Nint(\mathcal{C})) = Nint(\mathcal{C})$.
(iii) To prove $N\alpha int(NS_\alpha int(\mathcal{C})) = NS_\alpha int(N\alpha int(\mathcal{C})) = N\alpha int(\mathcal{C})$.
Since $N\alpha int(\mathcal{C})$ is $N\alpha$-O.S., therefore $N\alpha int(\mathcal{C})$ is $NS_\alpha$-O.S.. Therefore by proposition (4.3):
$N\alpha int(\mathcal{C}) = NS_\alpha int(N\alpha int(\mathcal{C}))$ ................................................................................ (1)
Now, to prove $N\alpha int(\mathcal{C}) = N\alpha int(NS_\alpha int(\mathcal{C}))$.
Since $N\alpha int(\mathcal{C}) \subseteq NS_\alpha int(\mathcal{C}) \Rightarrow N\alpha int(N\alpha int(\mathcal{C})) \subseteq N\alpha int(NS_\alpha int(\mathcal{C})) \Rightarrow N\alpha int(\mathcal{C}) \subseteq N\alpha int(NS_\alpha int(\mathcal{C}))$.
Also, $NS_\alpha int(\mathcal{C}) \subseteq \mathcal{C} \Rightarrow N\alpha int(NS_\alpha int(\mathcal{C})) \subseteq N\alpha int(\mathcal{C})$. Hence:
$N\alpha int(\mathcal{C}) = N\alpha int(NS_\alpha int(\mathcal{C}))$ ................................................................................ (2)
Therefore by (1) and (2), we get $N\alpha int(NS_\alpha int(\mathcal{C})) = NS_\alpha int(N\alpha int(\mathcal{C})) = N\alpha int(\mathcal{C})$.
(iv) To prove $Ncl(NS_\alpha cl(\mathcal{C})) = NS_\alpha cl(Ncl(\mathcal{C})) = Ncl(\mathcal{C})$.
We know that $Ncl(\mathcal{C})$ is a $N$-C.S., so it is $NS_\alpha$-C.S.. Hence by proposition (4.3), we have:
$Ncl(\mathcal{C}) = NS_\alpha cl(Ncl(\mathcal{C}))$ ........................................................................................ (1)

To prove $Ncl(\mathcal{C}) = Ncl(NS_\alpha cl(\mathcal{C}))$. Since $NS_\alpha cl(\mathcal{C}) \subseteq Ncl(\mathcal{C})$ (by part (i)).
Then $Ncl(NS_\alpha cl(\mathcal{C})) \subseteq Ncl(Ncl(\mathcal{C})) = Ncl(\mathcal{C}) \Rightarrow Ncl(NS_\alpha cl(\mathcal{C})) \subseteq Ncl(\mathcal{C})$.
Since $\mathcal{C} \subseteq NS_\alpha cl(\mathcal{C}) \subseteq Ncl(NS_\alpha cl(\mathcal{C}))$, then $\mathcal{C} \subseteq Ncl(NS_\alpha cl(\mathcal{C}))$. Hence $Ncl(\mathcal{C}) \subseteq Ncl(Ncl(NS_\alpha cl(\mathcal{C}))) = Ncl(NS_\alpha cl(\mathcal{C})) \Rightarrow Ncl(\mathcal{C}) \subseteq Ncl(NS_\alpha cl(\mathcal{C}))$ and therefore:
$Ncl(\mathcal{C}) = Ncl(NS_\alpha cl(\mathcal{C}))$..............................................................................(2)
Now, by (1) and (2), we get that $Ncl(NS_\alpha cl(\mathcal{C})) = NS_\alpha cl(Ncl(\mathcal{C}))$.
Hence $Ncl(NS_\alpha cl(\mathcal{C})) = NS_\alpha cl(Ncl(\mathcal{C})) = Ncl(\mathcal{C})$.
(vii) To prove $NS_\alpha int(\mathcal{C}) = \mathcal{C} \cap Ncl(Nint(Ncl(Nint(\mathcal{C}))))$.
Since $NS_\alpha int(\mathcal{C}) \in NS_\alpha O(\mathcal{U}, \mathcal{M}) \Rightarrow NS_\alpha int(\mathcal{C}) \subseteq Ncl(Nint(Ncl(Nint(NS_\alpha int(\mathcal{C})))))$
$= Ncl(Nint(Ncl(Nint(\mathcal{C}))))$ (by part (ii)).
Hence $NS_\alpha int(\mathcal{C}) \subseteq Ncl(Nint(Ncl(Nint(\mathcal{C}))))$, also $NS_\alpha int(\mathcal{C}) \subseteq \mathcal{C}$. Then:
$NS_\alpha int(\mathcal{C}) \subseteq \mathcal{C} \cap Ncl(Nint(Ncl(Nint(\mathcal{C}))))$ ............................................................(1)
To prove $\mathcal{C} \cap Ncl(Nint(Ncl(Nint(\mathcal{C}))))$ is a $NS_\alpha$-O.S. contained in $\mathcal{C}$.
It is clear that $\mathcal{C} \cap Ncl(Nint(Ncl(Nint(\mathcal{C})))) \subseteq Ncl(Nint(Ncl(Nint(\mathcal{C}))))$ and also it is clear that $Nint(\mathcal{C}) \subseteq Ncl(Nint(\mathcal{C})) \Rightarrow Nint(Nint(\mathcal{C})) \subseteq Nint(Ncl(Nint(\mathcal{C})))$
$\Rightarrow Nint(\mathcal{C}) \subseteq Nint(Ncl(Nint(\mathcal{C}))) \Rightarrow Ncl(Nint(\mathcal{C})) \subseteq Ncl(Nint(Ncl(Nint(\mathcal{C}))))$ and $Nint(\mathcal{C}) \subseteq Ncl(Nint(\mathcal{C})) \Rightarrow Nint(\mathcal{C}) \subseteq Ncl(Nint(Ncl(Nint(\mathcal{C}))))$ and $Nint(\mathcal{C}) \subseteq \mathcal{C}$
$\Rightarrow Nint(\mathcal{C}) \subseteq \mathcal{C} \cap Ncl(Nint(Ncl(Nint(\mathcal{C}))))$.
We get $Nint(\mathcal{C}) \subseteq \mathcal{C} \cap Ncl(Nint(Ncl(Nint(\mathcal{C})))) \subseteq Ncl(Nint(Ncl(Nint(\mathcal{C}))))$.
Hence $\mathcal{C} \cap Ncl(Nint(Ncl(Nint(\mathcal{C}))))$ is a $NS_\alpha$-O.S. (by proposition (4.3)).
Also, $\mathcal{C} \cap Ncl(Nint(Ncl(Nint(\mathcal{C}))))$ is contained in $\mathcal{C}$. Then $\mathcal{C} \cap Ncl(Nint(Ncl(Nint(\mathcal{C}))))$
$\subseteq NS_\alpha int(\mathcal{C})$ (since $NS_\alpha int(\mathcal{C})$ is the largest $NS_\alpha$-O.S. contained in $\mathcal{C}$). Hence:
$\mathcal{C} \cap Ncl(Nint(Ncl(Nint(\mathcal{C})))) \subseteq NS_\alpha int(\mathcal{C})$ ...................................................... (2)
By (1) and (2), $NS_\alpha int(\mathcal{C}) = \mathcal{C} \cap Ncl(Nint(Ncl(Nint(\mathcal{C}))))$.
(viii) To prove that $Nint(Ncl(\mathcal{C})) \subseteq NS_\alpha int(NS_\alpha cl(\mathcal{C}))$. Since $NS_\alpha cl(\mathcal{C})$ is a $NS_\alpha$-C.S., therefore $Nint(Ncl(Nint(Ncl(NS_\alpha cl(\mathcal{C}))))) \subseteq NS_\alpha cl(\mathcal{C})$ (by corollary (3.12)).
Hence $Nint(Ncl(\mathcal{C})) \subseteq Nint(Ncl(Nint(Ncl(\mathcal{C})))) \subseteq NS_\alpha cl(\mathcal{C})$ (by part (iv)).
Therefore, $NS_\alpha int(Nint(Ncl(\mathcal{C}))) \subseteq NS_\alpha int(NS_\alpha cl(\mathcal{C})) \Rightarrow$
$Nint(Ncl(\mathcal{C})) \subseteq NS_\alpha int(NS_\alpha cl(\mathcal{C}))$ (by part (ii)).

**Theorem 4.10:**

For any subset $\mathcal{C}$ of a N.T.S. $(\mathcal{U}, \tau_\mathcal{R}(\mathcal{M}))$. The following properties are equivalent:
(i) $\mathcal{C} \in NS_\alpha O(\mathcal{U}, \mathcal{M})$.
(ii) $\mathcal{K} \subseteq \mathcal{C} \subseteq Ncl(Nint(Ncl(\mathcal{K})))$, for some $N$-O.S. $\mathcal{K}$.
(iii) $\mathcal{K} \subseteq \mathcal{C} \subseteq Nsint(Ncl(\mathcal{K}))$, for some $N$-O.S. $\mathcal{K}$.
(iv) $\mathcal{C} \subseteq Nsint(Ncl(Nint(\mathcal{C})))$.

*Proof:*
$(i) \Rightarrow (ii)$ Let $\mathcal{C} \in NS_\alpha O(\mathcal{U}, \mathcal{M})$, then $\mathcal{C} \subseteq Ncl(Nint(Ncl(Nint(\mathcal{C}))))$ and $Nint(\mathcal{C}) \subseteq \mathcal{C}$.
Hence $\mathcal{K} \subseteq \mathcal{C} \subseteq Ncl(Nint(Ncl(\mathcal{K})))$, where $\mathcal{K} = Nint(\mathcal{C})$.
$(ii) \Rightarrow (iii)$ Suppose $\mathcal{K} \subseteq \mathcal{C} \subseteq Ncl(Nint(Ncl(\mathcal{K})))$, for some $N$-O.S. $\mathcal{K}$.
But $Nsint(Ncl(\mathcal{K})) = Ncl(Nint(Ncl(\mathcal{K})))$ (by lemma (2.6)).
Then $\mathcal{K} \subseteq \mathcal{C} \subseteq Nsint(Ncl(\mathcal{K}))$, for some $N$-O.S. $\mathcal{K}$.
$(iii) \Rightarrow (iv)$ Suppose that $\mathcal{K} \subseteq \mathcal{C} \subseteq Nsint(Ncl(\mathcal{K}))$, for some $N$-O.S. $\mathcal{K}$.
Since $\mathcal{K}$ is a $N$-O.S. contained in $\mathcal{C}$. Then $\mathcal{K} \subseteq Nint(\mathcal{C}) \Rightarrow Ncl(\mathcal{K}) \subseteq Ncl(Nint(\mathcal{C}))$
$\Rightarrow Nsint(Ncl(\mathcal{K})) \subseteq Nsint(Ncl(Nint(\mathcal{C})))$. But $\mathcal{C} \subseteq Nsint(Ncl(\mathcal{K}))$ (by hypothesis),

then $C \subseteq Nsint(Ncl(Nint(C)))$.
$(iv) \Rightarrow (i)$ Let $C \subseteq Nsint(Ncl(Nint(C)))$.
But $Nsint(Ncl(Nint(C))) = Ncl(Nint(Ncl(Nint(C))))$ (by lemma (2.6)).
Hence $C \subseteq Ncl(Nint(Ncl(Nint(C)))) \Rightarrow C \in NS_\alpha O(\mathcal{U}, \mathcal{M})$.

**Corollary 4.11:**

For any subset $\mathcal{D}$ of a N. T. S. $(\mathcal{U}, \tau_\mathcal{R}(\mathcal{M}))$, the following properties are equivalent:
(i) $\mathcal{D} \in NS_\alpha C(\mathcal{U}, \mathcal{M})$.
(ii) $Nint(Ncl(Nint(\mathcal{F}))) \subseteq \mathcal{D} \subseteq \mathcal{F}$, for some $\mathcal{F}$ $N$-C. S..
(iii) $Nscl(Nint(\mathcal{F})) \subseteq \mathcal{D} \subseteq \mathcal{F}$, for some $\mathcal{F}$ $N$-C. S..
(iv) $Nscl(Nint(Ncl(\mathcal{D}))) \subseteq \mathcal{D}$.

*Proof:*
$(i) \Rightarrow (ii)$ Let $\mathcal{D} \in NS_\alpha C(\mathcal{U}, \mathcal{M}) \Rightarrow Nint(Ncl(Nint(Ncl(\mathcal{D})))) \subseteq \mathcal{D}$ (by corollary (3.12))
and $\mathcal{D} \subseteq Ncl(\mathcal{D})$. Hence we get $Nint(Ncl(Nint(Ncl(\mathcal{D})))) \subseteq \mathcal{D} \subseteq Ncl(\mathcal{D})$.
Therefore $Nint(Ncl(Nint(\mathcal{F}))) \subseteq \mathcal{D} \subseteq \mathcal{F}$, where $\mathcal{F} = Ncl(\mathcal{D})$.
$(ii) \Rightarrow (iii)$ Let $Nint(Ncl(Nint(\mathcal{F}))) \subseteq \mathcal{D} \subseteq \mathcal{F}$, for some $\mathcal{F}$ $N$-C. S..
But $Nint(Ncl(Nint(\mathcal{F}))) = Nscl(Nint(\mathcal{F}))$ (by lemma (2.6)).
Hence $Nscl(Nint(\mathcal{F})) \subseteq \mathcal{D} \subseteq \mathcal{F}$, for some $\mathcal{F}$ $N$-C. S..
$(iii) \Rightarrow (iv)$ Let $Nscl(Nint(\mathcal{F})) \subseteq \mathcal{D} \subseteq \mathcal{F}$, for some $\mathcal{F}$ $N$-C. S..
Since $\mathcal{D} \subseteq \mathcal{F}$ (by hypothesis), hence $Ncl(\mathcal{D}) \subseteq \mathcal{F} \Rightarrow Nint(Ncl(\mathcal{D}) \subseteq Nint(\mathcal{F}) \Rightarrow Nscl(Nint(Ncl(\mathcal{D}))) \subseteq Nscl(Nint(\mathcal{F})) \subseteq \mathcal{D} \Rightarrow Nscl(Nint(Ncl(\mathcal{D}))) \subseteq \mathcal{D}$.
$(iv) \Rightarrow (i)$ Let $Nscl(Nint(Ncl(\mathcal{D}))) \subseteq \mathcal{D}$.
But $Nscl(Nint(Ncl(\mathcal{D}))) = Nint(Ncl(Nint(Ncl(\mathcal{D}))))$ (by lemma (2.6)).
Hence $Nint(Ncl(Nint(Ncl(\mathcal{D})))) \subseteq \mathcal{D} \Rightarrow \mathcal{D} \in NS_\alpha C(\mathcal{U}, \mathcal{M})$.

## 5. CONCLUSION

The class of $NS_\alpha$-O. S. defined using $N\alpha$-O. S. forms a nano topology and lay between the class of $N$-O. S. and the class of $Ns$-O. S.. The $NS_\alpha$-O. S. can be used to derive a new decomposition of nano continuity, nano compactness, and nano connectedness.